\documentclass{amsart}
\usepackage{amsrefs,mathrsfs,math}
\newcommand\K{{\mathbb K}}
\numberwithin{equation}{section}
\newcommand\mcolonn{m:n}

\begin{document}
\title{A converse to Moore's theorem on cellular automata}
\author{Laurent Bartholdi}
\date{typeset \today; last timestamp 20070922}
\address{\'Ecole Polytechnique F\'ed\'erale de Lausanne (EPFL),
  Institut de Math\'ematiques B (IMB), 1015 Lausanne, Switzerland}
\email{laurent.bartholdi@gmail.com}
\thanks{This work was partially supported by a CNRS visiting position
  at Universit\'e de Provence, Marseille}
\begin{abstract}
  We prove a converse to Moore's ``Garden-of-Eden'' theorem: a group
  $G$ is amenable if and only if all cellular automata living on $G$
  that admit mutually erasable patterns also admit gardens of Eden.

  It had already been conjectured
  in~\cite{ceccherini-m-s:ca}*{Conjecture 6.2} that amenability could
  be characterized by cellular automata. We prove the first part of
  that conjecture.
\end{abstract}
\maketitle

\section{Introduction}
\begin{defn}\label{def:ca}
  Let $G$ be a group.  A finite \emph{cellular automaton} on $G$ is a
  map $\theta:Q^S\to Q$, where $Q$, the \emph{state set}, is a finite
  set, and $S$ is a finite subset of $G$.
\end{defn}
Note that usually $G$ is infinite; much of the theory holds trivially
if $G$ is finite. $S$ could be taken to be a generating set of $G$,
though this is not a necessity.

A cellular automaton should be thought of as a highly regular animal,
composed of many cells labeled by $G$, each in a state $\in Q$. Each
cell ``sees'' its neighbours as defined by $S$, and ``evolves''
according to its neighbours' states.

More formally: a \emph{configuration} is a map $\phi:G\to Q$. The
\emph{evolution} of the automaton $\theta:Q^S\to Q$ is the self-map
$\Theta:Q^G\to Q^G$ on configurations, defined by
\[\Theta(\phi)(x)=\theta(s\mapsto\phi(xs)).\]

Two properties of cellular automata received special attention. Let us
call \emph{patch} the restriction of a configuration to a finite
subset $Y\subseteq G$. On the one hand, there can exist patches that
never appear in the image of $\Theta$. These are called \emph{Garden
  of Eden} (GOE), the biblical metaphor expressing the notion of
paradise lost forever.

On the other hand, $\Theta$ can be non-injective in a strong sense:
there can exist patches $\phi'_1\neq\phi'_2\in Q^Y$ such that, however
one extends $\phi'_1$ to a configuration $\phi_1$, if one extends
$\phi'_2$ similarly (i.e.\ in such a way that $\phi_1$ and $\phi_2$
have the same restriction to $G\setminus Y$) then
$\Theta(\phi_1)=\Theta(\phi_2)$. These patches $\phi'_1,\phi'_2$ are
called \emph{Mutually Erasable Patterns} (MEP).
Equivalently\footnote{In the non-trivial direction, let
  $\phi_1,\phi_2$ differ on a non-empty finite set $F$; set $Y=F(S\cup
  S^{-1})$ and let $\phi'_1,\phi'_2$ be the restriction of
  $\phi_1,\phi_2$ to $Y$ respectively.} there are two configurations
$\phi_1,\phi_2$ which differ on a non-empty finite set, with
$\Theta(\phi_1)=\Theta(\phi_2)$. The absence of MEP is sometimes
called \emph{pre-injectivity}.

Cellular automata were initially considered on $G=\Z^n$. Celebrated
theorems by Moore and Myhill~\cites{moore:ca,myhill:ca} prove that, in
this context, a cellular automaton admits GOE if and only if it admits
MEP. This result was generalized by Mach\`\i\ and
Mignosi~\cite{machi-m:ca} to $G$ of subexponential growth, and by
Ceccherini, Mach\`\i\ and Scarabotti~\cite{ceccherini-m-s:ca} to $G$
amenable.

We prove that this last result is essentially optimal, and yields a
characterization of amenable groups:
\begin{thm}\label{thm:main}
  Let $G$ be a group. Then the following are equivalent:
  \begin{enumerate}
  \item the group $G$ is amenable;\label{thm:1}
  \item all cellular automata on $G$ that admit MEP also admit GOE.\label{thm:2}
  \end{enumerate}
\end{thm}
Schupp had already asked in~\cite{schupp:arrays}*{Question~1} in which
precise class of groups the Moore-Myhill theorem holds.

Ceccherini et al.\ write in~\cite{ceccherini-m-s:ca}\footnote{I changed
slightly their wording to match this paper's}:
\begin{conj}[\cite{ceccherini-m-s:ca}*{Conjecture 6.2}]\label{conj:cms}
  Let $G$ be a non-amenable finitely generated group. Then for any
  finite and symmetric generating set $S$ for $G$ there exist cellular
  automata $\theta_1,\theta_2$ \emph{with that $S$} such that
  \begin{itemize}
  \item In $\theta_1$ there are MEP but no GOE;
  \item In $\theta_2$ there are GOE but no MEP.
  \end{itemize}
\end{conj}

As a first step, we will prove Theorem~\ref{thm:main}, in which we
allow ourselves to choose an appropriate subset $S$ of $G$. Next, we
extend a little the construction to answer the first part of
Conjecture~\ref{conj:cms}:
\begin{thm}\label{thm:cms}
  Let $G=\langle S\rangle$ be a finitely generated, non-amenable
  group. Then there exists a cellular automaton $\theta:Q^S\to Q$ that
  has MEP but no GOE.
\end{thm}

We conclude that the property of ``satisfying Moore's theorem'' is
independent of the generating set, a fact which was not obvious
\emph{a priori}.

\section{Proof of Theorem~\ref{thm:main}}
The implication $\eqref{thm:1}\Rightarrow\eqref{thm:2}$ has been
proven by Ceccherini et al.; see
also~\cite{gromov:endomorphisms}*{\S8} for a slicker proof. We prove
the converse.

Let us therefore be given a non-amenable group $G$. Let us also, as a
first step, be given a large enough finite subset $S$ of $G$.  Then
there exists a ``bounded propagation $2:1$ compressing vector field''
on $G$: a map $f:G\to G$ such that $f(x)^{-1}x\in S$ and
$\#f^{-1}(x)=2$ for all $x\in G$.

We construct the following automaton $\theta$. Its stateset is
\[Q=S\times\{0,1\}\times S.
\]
Order $S$ in an arbitrary manner, and choose an arbitrary $q_0\in Q$.
Define $\theta:Q^S\to Q$ as follows:
\begin{equation}\label{eq:theta}
  \theta(\phi)=\begin{cases}
    (p,\alpha,q) & \text{for the minimal pair $s<t$ in $S$ with}
      \left\{\!\begin{array}{l}      \phi(s)=(s,\alpha,p),\\\phi(t)=(t,\beta,q),
      \end{array}\right.\\
    q_0 & \text{if no such $s,t$ exist}.
  \end{cases}
\end{equation}

\subsection{$\Theta$ is surjective}
Namely, $\theta$ does not admit GOE. Let indeed $\phi$ be any
configuration. We construct a configuration $\psi$ with
$\Theta(\psi)=\phi$.

Consider in turn all $x\in G$; write $\phi(x)=(p,\alpha,q)$, and
$f^{-1}(x)=\{xs,xt\}$ for some $s,t\in S$ ordered as $s<t$. Set then
\begin{equation}\label{eq:psi}
  \psi(xs)=(s,\alpha,p),\qquad\psi(xt)=(t,0,q).
\end{equation}

Note that $\psi(z)=(f^{-1}(z)z,*,*)$ for all $z\in G$. Since
$\#f^{-1}(z)=2$ for all $z\in G$, it is clear that, for every $x\in
G$, there are exactly two $s\in S$ such that $\psi(xs)=(s,*,*)$; call
them $s,t$, ordered such that $\psi(xs)=(s,\alpha,p)$ and
$\psi(xt)=(t,0,q)$.  Then $\Theta(\psi)(x)=(p,\alpha,q)$, so
$\Theta(\psi)=\phi$.

\subsection{$\Theta$ is not pre-injective}
Namely, $\theta$ admits MEP. Let indeed $\phi:G\to Q$ be any
configuration; then construct $\psi$ following~\eqref{eq:psi}, and
define $\psi'$ as follows. Choose any $y\in G$, write
$\phi(y)=(p,\alpha,q)$, and write $f^{-1}(y)=\{ys,yt\}$ for some
$s,t\in S$, ordered as $s<t$.  Define
$\psi':G\to Q$ by
\[\psi'(x)=\begin{cases}
  \psi(x) & \text{ if }x\neq yt,\\
  (t,1,q) & \text{ if }x=yt.
\end{cases}
\]
Then $\psi$ and $\psi'$ differ only at $yt$; and
$\Theta(\psi)=\Theta(\psi')$ because the value of $\beta$ is unused
in~\eqref{eq:theta}. We conclude that $\theta$ has MEP.

\section{Proof of Theorem~\ref{thm:cms}}
We begin by a new formulation of amenability for finitely generated
groups:
\begin{lem}\label{lem:mn}
  Let $G$ be a finitely generated group. The following are equivalent:
  \begin{enumerate}
  \item the group $G$ is not amenable;
  \item for every generating set $S$ of $G$, there exist $m>n\in\N$
    and a ``$\mcolonn$ compressing correspondence on $G$ with
    propagation $S$''; i.e.\ a function $f:G\times G\to\N$ such that
    \begin{align}
      \forall y\in G:&\quad\sum_{x\in G}f(x,y)=m,\\
      \forall x\in G:&\quad\sum_{y\in G}f(x,y)=n,\\
      \forall x,y\in G:&\quad f(x,y)\neq0\Rightarrow y\in xS.
    \end{align}
  \end{enumerate}
\end{lem}
Note that this definition generalizes the notion of ``$2:1$
compressing vector field'' introduced above.

\begin{proof}
  For the forward direction, assuming that $G$ is non-amenable, there
  exists a rational $m/n>1$ such that every finite $F\subseteq G$
  satisfies
  \[\#(FS)\ge m/n\#F.\] Construct the following bipartite oriented
  graph: its vertex set is $G\times\{1,\dots,m\}\sqcup
  G\times\{-1,\dots,-n\}$. There is an edge from $(g,i)$ to $(gs,-j)$
  for all $s\in S$ and all $i\in\{1,\dots,m\},j\in\{1,\dots,m\}$. By
  hypothesis, this graph satisfies: every finite $F\subseteq
  G\times\{1,\dots,m\}$ has at least $\#F$ neighbours. Since $m>n$ and
  multiplication by a generator is a bijection, every finite
  $F\subseteq G\times\{-1,\dots,-n\}$ also has at least $\#F$
  neighbours.

  We now invoke the Hall-Rado theorem~\cite{mirsky:transversal}: if a
  bipartite graph is such that every subset of any of the parts has as
  many neighbours as its cardinality, then there exists a ``perfect
  matching'' --- a subset $I$ of the edge set of the graph such that
  every vertex is contained in precisely one edge in $I$. Set then
  \[f(x,y)=\#\{(i,j)\in\{1,\dots,m\}\times\{1,\dots,n\}:\;I\text{ contains the edge from }(x,i)\text{ to }(y,-j)\}.
  \]

  For the backward direction: if $G$ is amenable, then there exists an
  invariant measure on $G$, hence on bounded natural-valued functions
  on $G$. Let $f$ be a bounded-propagation $\mcolonn$ compressing
  correspondence; then
  \[m=\sum_{x\in G}\int_{\{x\}\times G}f=\sum_{y\in
    G}\int_{G\times\{y\}}f=n,
  \]
  contradicting $m>n$.
\end{proof}

Let now $G=\langle S\rangle$ be a non-amenable group, and apply
Lemma~\ref{lem:mn} to $G=\langle S^{-1}\rangle$, yielding $m>n\in\N$
and a contracting $\mcolonn$ correspondence $f$. Consider the
following cellular automaton $\theta$, with stateset
\[Q=(S\times\{0,1\}\times S^n)^n.
\]
Choose $q_0\in Q$, and give a total ordering to
$S\times\{1,\dots,n\}$.

Consider $\phi\in Q^S$. To define $\theta(\phi)$, let
$(s_1,k_1)<\dots<(s_m,k_m)$ be the lexicographically minimal sequence
in $(S\times\{1,\dots,n\})^m$ such that
\[\phi(s_j)_{k_j}=(s_j,\alpha_j,t_{j,1},\dots,t_{j,n})\in
S\times\{0,1\}\times S^n\quad\text{ for }j=1,\dots,m.
\]
If no such $s_1,k_1,\dots,s_m,k_m$ exist, set $\theta(\phi)=q_0$;
otherwise, set
\begin{equation}\label{eq:theta2}
\theta(\phi)=((t_{1,1},\alpha_1,t_{2,1},\dots,t_{n+1,1}),\dots,(t_{1,n},\alpha_n,t_{2,n},\dots,t_{n+1,n}))\in Q.
\end{equation}

The same arguments as before apply. Given $\phi:G\to Q$, we construct
$\psi:G\to Q$ such that $\Theta(\psi)=\phi$, as follows. We think of
the co\"ordinates $\psi(x)_k$ of $\psi(x)$ as $n$ ``slots'', initially
all ``free''. By definition, $\#f^{-1}(x)=m$ for all $x\in G$, while
$\#f(x)=n$. Consider in turn all $x\in G$; write
$f^{-1}(x)=\{xs_1,\dots,xs_m\}$, and let
$k_1,\dots,k_m\in\{1,\dots,n\}$ be ``free'' slots in
$\psi(xs_1),\dots,\psi(xs_m)$ respectively. By the definition of $f$,
there always exist sufficiently many free slots.

Mark now these slots as ``occupied''. Reorder $s_1,k_1,\dots,s_m,k_m$
in such a way that $(s_1,k_1,\dots,s_m,k_m)$ is minimal among its $m!$
permutations. Set then
\[\psi(xs_j)_{k_j}=(s_j,\alpha_j,t_{j,1},\dots,t_{j,n})\quad\text{ for }j=1,\dots,m,
\]
where $\alpha_{n+1},\dots,\alpha_m$ are taken to be arbitrary values
(say $0$ for definiteness) and
\[\phi(x)=((t_{1,1},\alpha_1,t_{2,1},\dots,t_{n+1,1}),\dots,(t_{1,n},\alpha_1,t_{2,n},\dots,t_{n+1,n})).\]
Finally, define $\psi$ arbitrarily on slots that are still ``free''.

It is clear that $\Theta(\psi)=\phi$, so $\theta$ does not have GOE.
On the other hand, $\theta$ has MEP as before, because the values of
$\alpha_j$ in~\eqref{eq:theta2} are not used for
$j\in\{n+1,\dots,m\}$.

\section{Remarks}
\subsection{$G$-sets} A cellular automaton could more generally be
defined on a right $G$-set $X$. There is a natural notion of
amenability for $G$-sets, but it is not clear exactly to which extent
Theorem~\ref{thm:main} can be generalized to that setting.

\subsection{Myhill's Theorem} It seems harder to produce
counterexamples to Myhill's theorem (``GOE imply MEP'') for arbitrary
non-amenable groups, although there exists an example on
$C=C_2*C_2*C_2$, due to Muller\footnote{University of Illinois 1976
  class notes}. Let us make our task even harder, and restrict
ourselves to linear automata over finite rings (so we assume $Q$ is a
module over a finite ring and the map $\theta:Q^S\to Q$ is linear).
The following approach seems promising.

\begin{conj}[Folklore? I learnt it from V.\ Guba]
  Let $G$ be a group. The following are equivalent:
  \begin{enumerate}
  \item The group $G$ is amenable;\\
  \item Let $\K$ be a field. Then $\K G$ admits right
    common multiples, i.e.\ for any $\alpha,\beta\in\K G$ there
    exist $\gamma,\delta\in\K G$ with
    $\alpha\gamma=\beta\delta$ and $(\gamma,\delta)\neq(0,0)$.
  \end{enumerate}
\end{conj}
The implication $(1)\Rightarrow(2)$ is easy, and follows from F\o
lner's criterion of amenability by linear algebra.

Assume now the ``hard'' direction of the conjecture. Given $G$
non-amenable, we may then find a finite field $\K$, and
$\alpha,\beta\in\K G$ that do not have a common right multiple.

Set $Q=\K^2$ with basis $(e_1,e_2)$, let $S$ contain the
inverses of the supports of $\alpha$ and $\beta$, and define the
cellular automaton $\theta:Q^S\to Q$ by
\[\theta(\phi)=\sum_{x\in
  G}\big(\alpha(x^{-1})\langle\phi(x)|e_1\rangle-\beta(x^{-1})\langle\phi(x)|e_2\rangle,0\big).
\]
Then $\theta$ has GOE, indeed any configuration not in $(\K\times0)^G$
is a GOE. On the other hand, if $\theta$ had MEP, then by linearity we
might as well assume $\Theta(\phi)=0$ for some non-zero
finitely-supported $\phi:G\to Q$. Write $\phi=(\gamma,\delta)$ in
co\"ordinates; then $\Theta(\phi)=0$ gives $\alpha\gamma=\beta\delta$,
showing that $\alpha,\beta$ actually did have a common right multiple.

Muller's example is in fact a special case of this construction, with
\[G=\langle x,y,z|x^2,y^2,z^2\rangle,
\]
$\K=\F[2]$, and $\alpha=x$, $\beta=y+z$.


\bib{ceccherini-m-s:ca}{article}{
  author={Ceccherini-Silberstein, T. G.},
  author={Mach{\`{\i }}, A.},
  author={Scarabotti, F.},
  title={Amenable groups and cellular automata},
  language={English, with English and French summaries},
  journal={Ann. Inst. Fourier (Grenoble)},
  volume={49},
  date={1999},
  number={2},
  pages={673--685},
  issn={0373-0956},
  review={\MR {1697376 (2000k:43001)}},
}

\bib{gromov:endomorphisms}{article}{
  author={Gromov, M.},
  title={Endomorphisms of symbolic algebraic varieties},
  journal={J. Eur. Math. Soc. (JEMS)},
  volume={1},
  date={1999},
  number={2},
  pages={109--197},
  issn={1435-9855},
  review={\MR {1694588 (2000f:14003)}},
}

\bib{machi-m:ca}{article}{
  author={Mach{\`{\i }}, Antonio},
  author={Mignosi, Filippo},
  title={Garden of Eden configurations for cellular automata on Cayley graphs of groups},
  journal={SIAM J. Discrete Math.},
  volume={6},
  date={1993},
  number={1},
  pages={44--56},
  issn={0895-4801},
  review={\MR {1201989 (95a:68084)}},
}

\bib{mirsky:transversal}{book}{
  author={Mirsky, Leonid},
  title={Transversal theory. An account of some aspects of combinatorial mathematics},
  series={Mathematics in Science and Engineering, Vol. 75},
  publisher={Academic Press},
  place={New York},
  date={1971},
  pages={ix+255},
  review={\MR {0282853 (44 \#87)}},
}

\bib{moore:ca}{inproceedings}{
  title={Machine models of self-reproduction},
  author={Moore, Edward F.},
  pages={17--33},
  booktitle={Essays on cellular automata},
  editor={Burks, Arthur W.},
  publisher={University of Illinois Press},
  place={Urbana, Ill.},
  date={1970},
  review={\MR {0299409 (45 \#8457)}},
}

\bib{myhill:ca}{article}{
  author={Myhill, John},
  title={The converse of Moore's Garden-of-Eden theorem},
  journal={Proc. Amer. Math. Soc.},
  volume={14},
  date={1963},
  pages={685--686},
  issn={0002-9939},
  review={\MR {0155764 (27 \#5698)}},
}

\bib{schupp:arrays}{article}{
  author={Schupp, Paul E.},
  title={Arrays, automata and groups---some interconnections},
  conference={ title={Automata networks}, address={Argel\`es-Village}, date={1986}, },
  book={ series={Lecture Notes in Comput. Sci.}, volume={316}, publisher={Springer}, place={Berlin}, },
  date={1988},
  pages={19--28},
  review={\MR {961274}},
}



\begin{bibsection}
\begin{biblist}
\bibselect{bartholdi,math}
\end{biblist}
\end{bibsection}
\end{document}